\documentclass[11pt]{article}%
\usepackage{amssymb}
\usepackage{amsfonts}
\usepackage{amsmath}
\usepackage{graphicx}%
\providecommand{\U}[1]{\protect\rule{.1in}{.1in}}
\newtheorem{theorem}{Theorem}

\begin{document}

\title{High-order fractional linear differential equations do not describe \\stochastic processes with finite state space}
\author{P. Van Mieghem\thanks{ Faculty of Electrical Engineering, Mathematics and
Computer Science, P.O Box 5031, 2600 GA Delft, The Netherlands; \emph{email}:
P.F.A.VanMieghem@tudelft.nl }}
\date{Delft University of Technology}
\maketitle

\begin{abstract}
The replacement of the classical derivative in the governing equation of a
continuous-time Markov process with finite state space by a fractional
derivative of arbitrary order $\alpha>1$ produces solutions that are
mathematically incompatible with the theory of stochastic processes.

\end{abstract}

\section{Introduction}

We focus on the extension of continuous-time Markov processes with finite
state space to cover a more general dependence of a stochastic process over
time. The essence of a continuous-time Markov process $X\left(  t\right)  $
with time $t\geq0$ lies in the semi-group property $P\left(  t+u\right)
=P\left(  t\right)  P\left(  u\right)  =P\left(  u\right)  P\left(  t\right)
$ of the transition probability matrix $P\left(  t\right)  $, with elements
$P_{jk}\left(  t\right)  =\Pr\left[  X\left(  t\right)  =j|X\left(  0\right)
=k\right]  $. The semi-group property results at any state in exponential
transition times, that feature the memoryless property, which induces the
one-step dependence of the event $\left\{  X\left(  t\right)  =k\right\}  $ on
$\left\{  X\left(  u\right)  =j\right\}  $ at the previous time $u<t$ at which
the Markov process changes. Real processes do not always possess the Markov
properties and motivate the search for a more realistic modeling.

Fractional calculus (see e.g. \cite{Podlubny1999},
\cite{Hilfer_Fractional_Calculus2000}, \cite{Sandev_Tomovski},
\cite{Gorenflo_2020}, \cite{Mainardi_bookFracCalculus2021}) may provide a
particular extension of continuous-time Markov processes with the main
attractive property of an exact, analytic solution
\cite{PVM_Mittag-Leffler_Gamma}. The fractional derivative defined in Section
\ref{sec_complex_fractional_derivative} is, however, far more complicated than
the classical derivative. While the properties of a fractional derivative are
mathematically studied and understood, its repercussions and its meaning to
physical phenomena, in particular stochastic processes, is far less understood.

Here, we highlight a serious limitation for stochastic processes: the
extension of the governing equation (\ref{differential_state_vector_s(t)}) of
a continuous-time Markovian process with $N$ states to a fractional setting is
shown to be mathematically consistent only for fractional order $0<\alpha
\leq1$. We demonstrate, generally, that the exact solution
(\ref{general_solution_s_a(t)_Re(a)<m}) of the general fractional order,
linear differential equation (\ref{fractional_differential_state_vector_s(t)})
for $m=\left\lceil \alpha\right\rceil >1$ cannot describe a stochastic process
with a finite state space. Previously in the literature (see e.g.
\cite{Mainardi_fractionalPoisson2004}, \cite{Mura_2008},
\cite{Beghin_Orsingher2009}, \cite{Georgiou_Kiss_ScalasPRE2015},
\cite{Georgiou_Scalas2022}), the confinement of the fractional order $\alpha$
to the interval $(0,1]$ was made based upon the complete monotonicity
\cite[Sec. 5.1]{PVM_Mittag-Leffler_Gamma} of the Mittag-Leffler function
$E_{a,b}\left(  z\right)  $ for $0<a\leq1$ and $b=1$. However, also for values
of $a$ outside the interval $(0,1]$, some combinations of Mittag-Leffler
functions are completely monotonous \cite[art. 91]{PVM_Mittag-Leffler_Gamma}.
Moreover, we deduce in \cite[art. 31]{PVM_Mittag-Leffler_Gamma} that the set
$\left\{  E_{a,b-j}\right\}  _{j\geq0}$ is independent, implying that a
function $f\left(  z\right)  $ can be written as a linear combination
$\sum_{j=0}^{\infty}c_{j}E_{a,b-j}\left(  z\right)  $ with the coefficients
$c_{j}$, depending on the function $f\left(  z\right)  $. Hence, real
non-negative functions such as probability density functions can be expressed
as a linear combination of Mittag-Leffler functions, while some of those
individual Mittag-Leffler functions might violate complete monotonicity. In
summary, complete monotonicity alone is insufficient to dismiss the
possibility that fractional linear differential equations of order $\alpha>1$
cannot describe stochastic processes.

Merely for readability and to motivate of the aim of this paper, we review
earlier established theory, such as the complex fractional derivative in
Section \ref{sec_complex_fractional_derivative} and the fractional extension
(\ref{fractional_differential_state_vector_s(t)}) of the Markov governing
equation (\ref{differential_state_vector_s(t)}) in Section
\ref{sec_fractional_continuous_time_MarkovProcesses}. Theorem
\ref{theorem_solution_frac_lin_diff} in Section
\ref{sec_proof_fractional_Markov_governing_equation} gives the general
solution of (\ref{fractional_differential_state_vector_s(t)}) and Section
\ref{sec_derivatives_law_tot_prob} formalizes the equation structure for
stochastic processes on a finite state space. Finally, our main result in
Section \ref{sec_formal_comparison} illustrates a matrix structure
incompatibility between a general-order fractional, linear process and a
continuous-time stochastic process on a finite state space. We conclude in
Section \ref{sec_conclusions} that the mere replacement of the classical
derivative by a fractional derivative needs a reconsideration.

\section{Fractional derivative of complex order}

\label{sec_complex_fractional_derivative}The $k$-th order derivative of a
function $f\left(  z\right)  $ can be expressed by Cauchy's integral
\cite{Titchmarshfunctions}
\begin{equation}
\left.  \frac{d^{k}f(z)}{dz^{k}}\right\vert _{z=z_{0}}=\frac{\Gamma\left(
1+k\right)  }{2\pi i}\int_{C(z_{0})}\frac{f(\omega)\;d\omega}{(\omega
-z_{0})^{k+1}} \label{k-th_derivative_f}%
\end{equation}
where $C(z_{0})$ is a contour around the point $z_{0}$ in a region of the
complex plane where the function $f\left(  z\right)  $ is analytic. In
\cite{PVM_fractional_derivative}, we have proposed and studied the
\emph{fractional derivative of complex order }$\alpha$, defined by the
\emph{formal} extension of the scope of the integer number $k$ in
(\ref{k-th_derivative_f}) to the complex number $\alpha$, as
\begin{equation}
D^{\alpha}f\left(  z_{0}\right)  \equiv\left.  \frac{d^{\alpha}f(z)}%
{dz^{\alpha}}\right\vert _{z=z_{0}}=\frac{\Gamma\left(  1+\alpha\right)
}{2\pi i}\int_{C(z_{0})}\frac{f(\omega)\;d\omega}{(\omega-z_{0})^{\alpha+1}}
\label{formal_extension_k_alpha_complex_derivative}%
\end{equation}
The fundamental problem in (\ref{formal_extension_k_alpha_complex_derivative})
with the extension of the integer $k$ to the complex $\alpha$ is generated by
the branch cut of the function $z^{\alpha}=e^{\alpha\log z}$ along the
negative real $z$-axis, causing that the function $\frac{f(\omega)}%
{(\omega-z_{0})^{\alpha+1}}$ is not analytic around the point $z_{0}$ and that
the contour $C\left(  z_{0}\right)  $ in
(\ref{formal_extension_k_alpha_complex_derivative}) cannot be closed over the
branch cut. Consequently, the formal extension
(\ref{formal_extension_k_alpha_complex_derivative}) is ill defined and its
reparations, as shown in \cite{PVM_fractional_derivative}, lead to infinitely
many possible definitions of the fractional derivative (see also
\cite[Appendix E]{Gorenflo_2020}). A general form for real $t$, that already
appeared in the early work of Riemann (see \cite[p. 6]%
{PVM_Mittag-Leffler_Gamma}), is%
\begin{equation}
D_{m,p}^{\alpha}f\left(  t\right)  =\left.  \frac{d^{\alpha}f(z)}{dz^{\alpha}%
}\right\vert _{z=t}=\frac{1}{\Gamma\left(  m-\alpha\right)  }\int_{p}^{t}%
\frac{f^{\left(  m\right)  }(x)}{\left(  t-x\right)  ^{\alpha+1-m}}%
dx\hspace{0.5cm}\text{for }\operatorname{Re}\left(  \alpha\right)  <m
\label{general_in_m_complex_fractional_derivative_no_jump}%
\end{equation}
where the real number $p<t$ and the $m$-th classical derivative $f^{\left(
m\right)  }(x)$ is assumed to be finite on the interval $x\in\left[
0,t\right]  $. The substitution $u=t-x$ in
(\ref{general_in_m_complex_fractional_derivative_no_jump}),%
\begin{equation}
\left.  \frac{d^{\alpha}f(z)}{dz^{\alpha}}\right\vert _{z=t}=\frac{1}%
{\Gamma\left(  m-\alpha\right)  }\int_{0}^{t-p}u^{m-\alpha-1}f^{\left(
m\right)  }(t-u)du\hspace{0.5cm}\text{for }\operatorname{Re}\left(
\alpha\right)  <m
\label{general_in_m_complex_fractional_derivative_no_jump_alternative}%
\end{equation}
further motivates to choose the point $p=0$. The simplest definition of the
fractional derivative is the case of $m=1$, real $\alpha$ and $p=0$ in
(\ref{general_in_m_complex_fractional_derivative_no_jump}),%
\begin{equation}
D_{0}^{\alpha}f\left(  t\right)  =\frac{1}{\Gamma\left(  1-\alpha\right)
}\int_{0}^{t}\frac{f^{\prime}(x)}{\left(  t-x\right)  ^{\alpha}}%
dx\hspace{0.5cm}\text{for }0<\alpha<1 \label{Caputo_fractional_derivative}%
\end{equation}
which is known as the case $m=1$ of the Caputo fractional derivative
\cite[Chapter E.2]{Gorenflo_2020}.

\section{Fractional extension of continuous-time Markov processes}

\label{sec_fractional_continuous_time_MarkovProcesses}We revisit the
generalization of continuous-time Markov processes \cite[Chapter
10]{PVM_PAComplexNetsCUP} to a \textquotedblleft fractional
setting\textquotedblright. The governing equation of a continuous-time
Markovian process with $N$ states is
\begin{equation}
\frac{ds\left(  t\right)  }{dt}=-Qs\left(  t\right)
\label{differential_state_vector_s(t)}%
\end{equation}
where $-Q$ is the infinitesimal generator in \cite[(10.13) after
transpose]{PVM_PAComplexNetsCUP}, equal to minus the $N\times N$ weighted
Laplacian matrix $Q$ of the underlying Markov graph, in which a node defines a
state in the Markov process. The $N\times1$ probability state
vector\footnote{In contrast to Markov theory \cite{PVM_PAComplexNetsCUP},
where vectors are usually written as $1\times N$ row vectors, we follow here
linear algebra as in \cite{PVM_graphspectra_second_edition} and write
$N\times1$ column vectors.} is
\[
s\left(  t\right)  =\left(  \Pr\left[  X\left(  t\right)  =1\right]
,\Pr\left[  X\left(  t\right)  =2\right]  ,\ldots,\Pr\left[  X\left(
t\right)  =N\right]  \right)
\]
and $s_{j}\left(  t\right)  =\Pr\left[  X\left(  t\right)  =j\right]  $
denotes the probability that the Markov process $X\left(  t\right)  $ is in
state $j$ at time $t$. The probability operator $\Pr\left[  .\right]  $ on
event $X$ measures the event $X$ by a non-negative, real number in the
interval $\left[  0,1\right]  $ so that $0\leq s_{j}\left(  t\right)  \leq1$.
Since the Markov process $X\left(  t\right)  $ at time $t$ must be in one of
the $N$ possible states, axiom 1 in probability theory \cite[p. 8]%
{PVM_PAComplexNetsCUP} states that $u^{T}s\left(  t\right)  =\sum_{j=1}^{N}%
\Pr\left[  X\left(  t\right)  =j\right]  =1$, where the unity vector
$u=\left(  1,1,\ldots,1\right)  $ is the all-one vector. Hence, it follows
from the governing equation (\ref{differential_state_vector_s(t)}) after
left-multiplication by $u^{T}$ that $u^{T}Q=0$, which is the characteristic
property of the weighted, generally asymmetric Laplacian $Q$. Indeed, the
unity or all-one vector $u$ is the left-eigenvector of any Laplacian matrix
$Q$ belonging to its zero eigenvalue.

We confine the analysis to an infinitesimal generator $Q$ that is fixed and
does not depend upon time $t$. In other words, the Markov \textquotedblleft
master\textquotedblright\ equation (\ref{differential_state_vector_s(t)}) is
equivalent to an $N$-th order linear differential equation with constant,
time-independent coefficients, whose general solution%
\[
s\left(  t\right)  =e^{-Qt}s\left(  0\right)
\]
possesses exponential functions at different eigenvalues of $Q$ as its
eigenfunctions (see \cite[p. 7-8]{PVM_Mittag-Leffler_Gamma}).

Instead of replacing the classical derivative in
(\ref{differential_state_vector_s(t)}) by the Caputo fractional derivative
(\ref{Caputo_fractional_derivative}) with limited order $0<\alpha<1$ as in
\cite{PVM_fractional_derivative}, we formally extend the governing equation
(\ref{differential_state_vector_s(t)}) with fixed Laplacian $Q$ to the general
fractional derivative
(\ref{general_in_m_complex_fractional_derivative_no_jump}), valid for
$\operatorname{Re}\left(  \alpha\right)  <m$,%
\begin{equation}
D_{m,0}^{\alpha}s_{\alpha}\left(  t\right)  =-Qs_{\alpha}\left(  t\right)
\label{fractional_differential_state_vector_s(t)}%
\end{equation}

\section{Solution of the general fractional linear differential equation
(\ref{fractional_differential_state_vector_s(t)})}

This section is devoted to the analytic solution of the $\alpha$-fractional
stochastic process $s_{\alpha}\left(  t\right)  $. We first derive in Section
\ref{sec_Laplace_transform_frac_diff_eq} a solution of
(\ref{fractional_differential_state_vector_s(t)}) under restricted conditions
by Laplace transformation. Inspired by that solution, we prove in Section
\ref{sec_proof_fractional_Markov_governing_equation} the general solution of
the general fractional linear differential equation
(\ref{fractional_differential_state_vector_s(t)}).

\subsection{Laplace transformation of
(\ref{fractional_differential_state_vector_s(t)})}

\label{sec_Laplace_transform_frac_diff_eq}Taking the Laplace transform of both
sides in (\ref{fractional_differential_state_vector_s(t)}) and using
(\ref{Laplace_transform_fractional_derivative_m_general}) yields, for
$\operatorname{Re}\left(  \alpha\right)  <m$,%
\begin{equation}
Q\mathcal{L}\left[  s_{\alpha}\left(  t\right)  \right]  =z^{\alpha-1}%
\sum_{l=0}^{m-1}z^{-l}s_{\alpha}^{\left(  l\right)  }\left(  0\right)
-z^{\alpha}\mathcal{L}\left[  s_{\alpha}(t)\right]
\label{Laplace_transform_fractional_Laplacian}%
\end{equation}
which is reorganized as%
\[
\left(  Q+z^{\alpha}I\right)  \mathcal{L}\left[  s_{\alpha}(t)\right]
=\sum_{l=0}^{m-1}z^{\alpha-1-l}s_{\alpha}^{\left(  l\right)  }\left(
0\right)
\]
Multiplying both sides by the matrix $\left(  Q+z^{\alpha}I\right)  ^{-1}$,
which we assume to exist, then results into the Laplace transform of the
probability state vector $s_{\alpha}(t)$,%
\begin{equation}
\mathcal{L}\left[  s_{\alpha}(t)\right]  =\sum_{l=0}^{m-1}z^{\alpha
-1-l}\left(  Q+z^{\alpha}I\right)  ^{-1}s_{\alpha}^{\left(  l\right)  }\left(
0\right)  \label{Laplace-tranformed_s_a_vector}%
\end{equation}

For simplicity\footnote{A general way invokes Jordan forms \cite{Meyer_matrix}%
, that lead to a considerably more involved derivation without adding new
insights. We have chosen to combine a Laplace approach -- a most common
approach in fractional calculus --, with a Taylor series approach. The latter
provides more clearly the restriction $m-1<\operatorname{Re}\left(
\alpha\right)  <m$ in Theorem \ref{theorem_solution_frac_lin_diff}.}, we first
proceed by assuming that the weighted Laplacian $Q$ is diagonalizable and the
eigenvalue decomposition \cite{PVM_graphspectra} is%
\[
Q=\sum_{k=1}^{N}\mu_{k}x_{k}y_{k}^{T}=XMY^{T}%
\]
where $\mu_{k}$, with $\operatorname{Re}\left(  \mu_{k}\right)  \geq0$, is the
eigenvalue belonging to the left-eigenvector $x_{k}$ and the right-eigenvector
$y_{k}$ of the weighted Laplacian $Q$ and we assume the ordering $\left\vert
\mu_{1}\right\vert \geq\left\vert \mu_{2}\right\vert \geq\cdots\geq\mu_{N}=0$.
The $N\times N$ matrix $X$ (and $Y$) contains the left-(right-) eigenvectors
in its columns and $M=$ diag$\left(  \mu\right)  $ is the $N\times N$ diagonal
matrix of the eigenvalue vector $\mu=\left(  \mu_{1},\mu_{2},\ldots,0\right)
$ with $k$-th diagonal element $M_{kk}=\mu_{k}$. Since $\sum_{k=1}^{N}%
x_{k}y_{k}^{T}=I$, because the right- and left-eigenvectors form a
bi-orthogonal set of vectors \cite[art. 237-239]%
{PVM_graphspectra_second_edition} and the corresponding right- and
left-eigenvector matrices obey $Y^{T}X=I$, we have%
\[
Q+z^{\alpha}I=\sum_{k=1}^{N}\left(  \mu_{k}+z^{\alpha}\right)  x_{k}y_{k}%
^{T}=X\left(  M+z^{\alpha}I\right)  Y^{T}%
\]
Provided that $\mu_{k}+z^{\alpha}\neq0$, we find%
\[
\left(  Q+z^{\alpha}I\right)  ^{-1}=\sum_{k=1}^{N}\frac{x_{k}y_{k}^{T}}%
{\mu_{k}+z^{\alpha}}=X\text{diag}\left(  \frac{1}{\mu+z^{\alpha}}\right)
Y^{T}%
\]
and the Laplace transform (\ref{Laplace-tranformed_s_a_vector}) becomes%
\begin{equation}
\mathcal{L}\left[  s_{\alpha}(t)\right]  =\sum_{k=1}^{N}x_{k}y_{k}^{T}%
\sum_{l=0}^{m-1}\frac{z^{\alpha-1-l}}{\mu_{k}+z^{\alpha}}s_{\alpha}^{\left(
l\right)  }\left(  0\right)  =X\sum_{l=0}^{m-1}\text{diag}\left(
\frac{z^{\alpha-1-l}}{\mu+z^{\alpha}}\right)  Y^{T}s_{\alpha}^{\left(
l\right)  }\left(  0\right)
\label{Laplace_transform_state_vector_eigendecomposition}%
\end{equation}
The Laplace transform (see e.g. \cite[art. 43]{PVM_Mittag-Leffler_Gamma}%
,\cite{Gorenflo_2020}) of the Mittag-Leffler function $E_{a,b}\left(
z\right)  =\sum_{k=0}^{\infty}\frac{z^{k}}{\Gamma\left(  b+ak\right)  }$,
which restricts $z$ by $\left\vert z^{a}\right\vert >\left\vert x\right\vert $
for $\operatorname{Re}\left(  a\right)  >0$ and $\operatorname{Re}\left(
b\right)  >0$,%
\begin{equation}
\int_{0}^{\infty}e^{-zt}t^{b-1}E_{a,b}\left(  xt^{a}\right)  dt=\frac{z^{a-b}%
}{z^{a}-x} \label{Laplace_transform_positive_argument}%
\end{equation}
indicates that the inverse Laplace transform is%
\[
\mathcal{L}^{-1}\left[  \frac{z^{\alpha-1-l}}{z^{\alpha}+\mu_{k}}\right]
=t^{l}E_{\alpha,l+1}\left(  -\mu_{k}t^{\alpha}\right)
\]
Inverse Laplace transformation of
(\ref{Laplace_transform_state_vector_eigendecomposition}) thus leads to%
\[
s_{\alpha}(t)=\sum_{l=0}^{m-1}t^{l}\left(  \sum_{k=1}^{N}E_{\alpha,l+1}\left(
-\mu_{k}t^{\alpha}\right)  x_{k}y_{k}^{T}\right)  s_{\alpha}^{\left(
l\right)  }\left(  0\right)
\]
Finally, after using $f\left(  Q\right)  =\sum_{k=1}^{N}f\left(  \mu
_{k}\right)  x_{k}y_{k}^{T}$ in \cite[(A.88)]{PVM_graphspectra_second_edition}
for any function $f$ defined on the eigenvalues of the matrix $Q$, we arrive
at the solution (\ref{general_solution_s_a(t)_Re(a)<m}) below of the
fractional matrix differential equation
(\ref{fractional_differential_state_vector_s(t)}) for $0<\operatorname{Re}%
\left(  \alpha\right)  <m$,which was already given without derivation in
\cite[Appendix C]{PVM_fractional_derivative}.

\subsection{General solution of
(\ref{fractional_differential_state_vector_s(t)}) for any $N\times N$ matrix
$Q$}

\label{sec_proof_fractional_Markov_governing_equation}The derivation in
Section \ref{sec_Laplace_transform_frac_diff_eq} has assumed that the matrix
$Q$ is diagonalizable. Here, we prove the general solution:

\begin{theorem}
\label{theorem_solution_frac_lin_diff}The general solution of the fractional,
linear differential equation $D_{m,0}^{\alpha}s_{\alpha}\left(  t\right)
=-Qs_{\alpha}\left(  t\right)  $ in
(\ref{fractional_differential_state_vector_s(t)}), where the matrix $Q$ can be
any $N\times N$ matrix with constant, time-independent elements and where $m$
is a non-negative integer, is%
\begin{equation}
s_{\alpha}(t)=\sum_{l=0}^{m-1}t^{l}E_{\alpha,l+1}\left(  -Qt^{\alpha}\right)
s_{\alpha}^{\left(  l\right)  }\left(  0\right)
\label{general_solution_s_a(t)_Re(a)<m}%
\end{equation}
with the restriction that the complex fractional order $\alpha$ must obey
$m-1<\operatorname{Re}\left(  \alpha\right)  <m$.
\end{theorem}

\textbf{Proof}: Substituting the suggestion $s_{\alpha}(t)=\sum_{l=0}%
^{m-1}t^{l}E_{\alpha,l+1}\left(  -Qt^{\alpha}\right)  s_{\alpha}^{\left(
l\right)  }\left(  0\right)  $ in (\ref{general_solution_s_a(t)_Re(a)<m}) into
the matrix differential equation $D_{m,0}^{\alpha}s_{\alpha}\left(  t\right)
=-Qs_{\alpha}\left(  t\right)  $ in
(\ref{fractional_differential_state_vector_s(t)}) with the explicit form of
the general fractional derivative
(\ref{general_in_m_complex_fractional_derivative_no_jump}) with $p=0$ first
needs the normal derivative $s_{\alpha}^{\left(  m\right)  }(x)=\frac{d^{m}%
}{dx^{m}}s_{\alpha}\left(  x\right)  $. The Taylor series $E_{a,b}\left(
z\right)  =\sum_{k=0}^{\infty}\frac{z^{k}}{\Gamma\left(  ak+b\right)  }$ of
the Mittag-Leffler function \cite{PVM_Mittag-Leffler_Gamma} converges for any
complex number $z$ for $\operatorname{Re}\left(  a\right)  >0$. A function of
a matrix is defined by a converging Taylor series \cite[Sec. 9.4; art.
232]{PVM_graphspectra_second_edition}. The Mittag-Leffler matrix function
$E_{a,b}\left(  Q\right)  =\sum_{k=0}^{\infty}\frac{Q^{k}}{\Gamma\left(
ak+b\right)  }$ for $\operatorname{Re}\left(  a\right)  >0$ and any complex
$b$ converges for any matrix $Q$.

After introducing the Taylor series in $s_{\alpha}(t)=\sum_{l=0}^{m-1}%
t^{l}E_{\alpha,l+1}\left(  -Qt^{\alpha}\right)  s_{\alpha}^{\left(  l\right)
}\left(  0\right)  $,%
\begin{equation}
s_{\alpha}(t)=\sum_{l=0}^{m-1}\sum_{k=0}^{\infty}\frac{\left(  -Q\right)
^{k}t^{\alpha k+l}}{\Gamma\left(  \alpha k+l+1\right)  }s_{\alpha}^{\left(
l\right)  }\left(  0\right)  \label{state_vector_ML_Taylor series}%
\end{equation}
the derivative $s_{\alpha}^{\left(  m\right)  }(x)=\frac{d^{m}}{dx^{m}%
}s_{\alpha}\left(  x\right)  $ is, with $\frac{d^{m}x^{\alpha k+l}}{dx^{m}%
}=\frac{\Gamma\left(  \alpha k+l\right)  }{\Gamma\left(  \alpha
k+l-m+1\right)  }x^{\alpha k+l-m}$, equal to%
\[
s_{\alpha}^{\left(  m\right)  }(x)=\sum_{l=0}^{m-1}\sum_{k=0}^{\infty}%
\frac{\left(  -Q\right)  ^{k}x^{\alpha k+l-m}}{\Gamma\left(  \alpha
k+l-m+1\right)  }s_{\alpha}^{\left(  l\right)  }\left(  0\right)
\]
While $s_{\alpha}(x)$ is regular for all real $x\geq0$, we observe that the
$m$-th derivative $s_{\alpha}^{\left(  m\right)  }(x)$ has a singularity at
$x=0$ if $\alpha k+l-m<0$ and $\alpha\notin\mathbb{N}$, because for
non-integer $\alpha k$, it holds that $\frac{1}{\Gamma\left(  \alpha
k+l-m+1\right)  }\neq0$. Thus, for $l\leq m-1$, the $k=0$ term disappears,
resulting in
\begin{equation}
s_{\alpha}^{\left(  m\right)  }(x)=\sum_{l=0}^{m-1}\sum_{k=1}^{\infty}%
\frac{\left(  -Q\right)  ^{k}x^{\alpha k+l-m}}{\Gamma\left(  \alpha
k+l-m+1\right)  }s_{\alpha}^{\left(  l\right)  }\left(  0\right)
\label{normal_derivative_state_vector_alpha}%
\end{equation}

The fractional derivative
(\ref{general_in_m_complex_fractional_derivative_no_jump}) thus becomes with
(\ref{normal_derivative_state_vector_alpha})
\begin{align*}
D_{m,0}^{\alpha}s_{\alpha}\left(  t\right)   &  =\frac{1}{\Gamma\left(
m-\alpha\right)  }\int_{0}^{t}\frac{s_{\alpha}^{\left(  m\right)  }%
(x)}{\left(  t-x\right)  ^{\alpha+1-m}}dx\\
&  =\frac{1}{\Gamma\left(  m-\alpha\right)  }\sum_{l=0}^{m-1}\sum
_{k=1}^{\infty}\frac{\left(  -Q\right)  ^{k}}{\Gamma\left(  \alpha
k+l-m+1\right)  }\left(  \int_{0}^{t}\frac{x^{\alpha k+l-m}dx}{\left(
t-x\right)  ^{\alpha+1-m}}\right)  s_{\alpha}^{\left(  l\right)  }\left(
0\right)
\end{align*}
The integral%
\[
\int_{0}^{t}\frac{x^{\alpha k+l-m}dx}{\left(  t-x\right)  ^{\alpha+1-m}%
}=t^{-\alpha-1+m}\int_{0}^{t}\left(  1-\frac{x}{t}\right)  ^{-\alpha
+m-1}x^{k\alpha+l-m}dx=t^{\left(  k-1\right)  \alpha+l}\int_{0}^{1}\left(
1-u\right)  ^{m-\alpha-1}u^{k\alpha+l-m}du
\]
equals, provided that $k\operatorname{Re}\left(  \alpha\right)  +l-m+1>0$,
\[
\int_{0}^{t}\frac{x^{k\alpha-1}}{\left(  t-x\right)  ^{\alpha}}dx=t^{\left(
k-1\right)  \alpha+l}\frac{\Gamma\left(  m-\alpha\right)  \Gamma\left(
k\alpha+l-m+1\right)  }{\Gamma\left(  \left(  k-1\right)  \alpha+l+1\right)  }%
\]
where we have used the Beta-function integral \cite[Chapter 17]%
{PVM_Mittag-Leffler_Gamma}, valid for $\operatorname{Re}\left(  z\right)  >0$
and $\operatorname{Re}\left(  q\right)  >0$,
\[
B\left(  z,q\right)  =\frac{\Gamma\left(  z\right)  \Gamma\left(  q\right)
}{\Gamma\left(  z+q\right)  }=\int_{0}^{1}u^{z-1}\left(  1-u\right)  ^{q-1}du
\]
Hence, the first validity condition of the Beta-integral requires that
$\operatorname{Re}\left(  \alpha\right)  >\frac{m-1-l}{k}$ for $k\geq1$, which
amounts to the requirement that $m-1<$ $\operatorname{Re}\left(
\alpha\right)  $ for all $k\geq1$, while the second condition of the
Beta-integral requires that $\operatorname{Re}\left(  \alpha\right)  <m$.
Thus, for $m-1<$ $\operatorname{Re}\left(  \alpha\right)  <m$, it holds that%
\begin{align*}
D_{m,0}^{\alpha}s_{\alpha}\left(  t\right)   &  =\frac{1}{\Gamma\left(
m-\alpha\right)  }\sum_{l=0}^{m-1}\sum_{k=1}^{\infty}\frac{\left(  -Q\right)
^{k}}{\Gamma\left(  \alpha k+l-m+1\right)  }\left(  t^{\left(  k-1\right)
\alpha+l}\frac{\Gamma\left(  m-\alpha\right)  \Gamma\left(  k\alpha
+l-m+1\right)  }{\Gamma\left(  \left(  k-1\right)  \alpha+l+1\right)
}\right)  s_{\alpha}^{\left(  l\right)  }\left(  0\right) \\
&  =\sum_{l=0}^{m-1}\sum_{k=1}^{\infty}\frac{\left(  -Q\right)  ^{k}t^{\left(
k-1\right)  \alpha+l}}{\Gamma\left(  \left(  k-1\right)  \alpha+l+1\right)
}s_{\alpha}^{\left(  l\right)  }\left(  0\right)  =-Q\sum_{l=0}^{m-1}%
\sum_{k=0}^{\infty}\frac{\left(  -Q\right)  ^{k}t^{k\alpha+l}}{\Gamma\left(
k\alpha+l+1\right)  }s_{\alpha}^{\left(  l\right)  }\left(  0\right)
\end{align*}
Finally with (\ref{state_vector_ML_Taylor series}), we arrive at
$D_{m,0}^{\alpha}s\left(  t\right)  =-Qs_{\alpha}(t)$, which proves Theorem
\ref{theorem_solution_frac_lin_diff}.$\hfill\square\medskip$

To some extent, the restriction $m-1<\operatorname{Re}\left(  \alpha\right)
<m$ in Theorem \ref{theorem_solution_frac_lin_diff} coincides with the fact
that a higher-order derivative in a physical governing equation describes a
different process. For example, the differential equations $\frac{ds}{dt}=-Qs$
and $\frac{d^{2}s}{dt^{2}}=-Qs$ have rather different solutions and appear in
different physical processes. A second fact is that the Mittag-Leffler
function $E_{a,b}\left(  z\right)  $, an entire function of order $\frac{1}%
{a}$, behaves seriously different in $a$ and different expressions (see
\cite[Chapter 4, e.g. art. 69]{PVM_Mittag-Leffler_Gamma}) also appear in the
regimes for $m-1<\operatorname{Re}\left(  \alpha\right)  <m$ for integer
$m=\left\lceil \operatorname{Re}\left(  \alpha\right)  \right\rceil $, where
$\left\lceil x\right\rceil $ is the smallest integer larger than or equal to
the real $x$.

\section{Stochastic processes}

\label{sec_D_law_total_prob}

\subsection{The law of total probability}

The law of total probability \cite{FellerI}, \cite[p. 23]%
{PVM_PAComplexNetsCUP}, applied to any continuous-time stochastic process
$Y\left(  t\right)  $ with a finite number of discrete states $n\in\left\{
1,2,\ldots,N\right\}  $ and time $t\geq0$, states that, for any real $v\geq0$,%
\begin{equation}
\Pr\left[  Y\left(  t+v\right)  =j\right]  =\sum_{n=1}^{N}\Pr\left[  \left.
Y\left(  t+v\right)  =j\right\vert Y\left(  v\right)  =n\right]  \Pr\left[
Y\left(  v\right)  =n\right]  \label{Prob_t+u_in_u_tot_law_prob}%
\end{equation}

Appendix \ref{sec_proof_law_total_prob} adds a proof of the law
(\ref{Prob_t+u_in_u_tot_law_prob}) of total probability.

\subsection{Derivatives of the law of total probability}

\label{sec_derivatives_law_tot_prob}The definition $\left(  \mathfrak{s}%
^{\left(  n\right)  }\left(  t\right)  \right)  _{j}=\frac{d^{n}}{dt^{n}}%
\Pr\left[  Y\left(  t\right)  =j\right]  $ and axiom 1 in probability theory
\cite[p. 8]{PVM_PAComplexNetsCUP}, $u^{T}\mathfrak{s}\left(  t\right)
=1=\sum_{n=1}^{N}\Pr\left[  Y\left(  v\right)  =n\right]  $ show, after
differentiation with respect to $t$, that any stochastic process obeys%
\begin{equation}
u^{T}\mathfrak{s}^{\left(  n\right)  }\left(  t\right)  =1_{\left\{
n=0\right\}  } \label{D_axiom_1_prob}%
\end{equation}
In other words, the vectors $\mathfrak{s}^{\left(  n\right)  }\left(
t\right)  $ for $n>0$ are orthogonal to the unity vector or all-one vector
$u=(1,1,\ldots1)$.

Axiom 1 in probability theory \textquotedblleft$\sum_{j=1}^{N}\Pr\left[
X_{\alpha}\left(  t\right)  =j\right]  =u^{T}s_{\alpha}\left(  t\right)  =1$
for all times $t$\textquotedblright, which equals (\ref{D_axiom_1_prob}) for
$n=0$, means that the process at time $t$ should be in one state in the
probability state space. Its application to $s_{\alpha}\left(  t\right)
=\sum_{l=0}^{m-1}t^{l}E_{\alpha,l+1}\left(  -Qt^{\alpha}\right)  s_{\alpha
}^{\left(  l\right)  }\left(  0\right)  $ in
(\ref{general_solution_s_a(t)_Re(a)<m}) is investigated in
\cite{PVM_fracDerivative_Markov_epidemics}, where we show that the solution
(\ref{general_solution_s_a(t)_Re(a)<m}) satisfies the necessary conditions
$u^{T}s_{\alpha}^{\left(  l\right)  }\left(  0\right)  =0$ for $0<l\leq m-1$.
However, we also show, by a counter example, that the necessary conditions are
not sufficient for $m=\left\lceil \alpha\right\rceil >1$ to guarantee that
$s_{\alpha}\left(  t\right)  $ is a probability state vector, i.e. that each
component $\left(  s_{\alpha}\left(  t\right)  \right)  _{n}$ for $1\leq n\leq
N$ is a probability that lies in the interval $\left[  0,1\right]  $. In the
sequel, we formalize the observation in
\cite{PVM_fracDerivative_Markov_epidemics}.

\section{Formal comparison between stochastic and fractional process}

\label{sec_formal_comparison}This Section \ref{sec_formal_comparison} compares
the characteristic properties of both a stochastic process on a finite state
space with a general-order fractional, linear process by recasting the
generalized state vector, defined in (\ref{def_sigma_vector}), in the
corresponding matrix representation (\ref{General_stochastic_eq}) and
(\ref{general_matrix_solution_general-order_frac_lin_diff_eq}). Section
\ref{sec_matrix_comparison} illustrates that the corresponding matrices
$\Upsilon(t,0)$ and $\digamma(t,0)$ exhibit a different structure, which
causes an incompatibility if a general-order fractional, linear process would
describe a stochastic process on a finite state space.

\subsection{Matrix equation of stochastic processes}

Leibniz' rule applied to the law (\ref{Prob_t+u_in_u_tot_law_prob}) of total
probability, expresses the $m$-th derivative as%
\begin{align*}
\frac{d^{m}}{dv^{m}}\Pr\left[  Y\left(  t+v\right)  =j\right]   &  =\sum
_{n=1}^{N}\frac{d^{m}}{dv^{m}}\left(  \Pr\left[  \left.  Y\left(  t+v\right)
=j\right\vert Y\left(  v\right)  =n\right]  \Pr\left[  Y\left(  v\right)
=n\right]  \right) \\
&  =\sum_{n=1}^{N}\sum_{l=0}^{m}\binom{m}{l}\frac{d^{m-l}}{dv^{m-l}}\Pr\left[
\left.  Y\left(  t+v\right)  =j\right\vert Y\left(  v\right)  =n\right]
\frac{d^{l}}{du^{l}}\Pr\left[  Y\left(  v\right)  =n\right]
\end{align*}
With the notation $\left(  \mathfrak{s}^{\left(  n\right)  }\left(  t\right)
\right)  _{j}=\frac{d^{n}}{dt^{n}}\Pr\left[  Y\left(  t\right)  =j\right]  $,
we obtain the $m$-th derivative of the law of total probability%
\begin{equation}
\left(  \mathfrak{s}^{\left(  m\right)  }\left(  t+v\right)  \right)
_{j}=\sum_{n=1}^{N}\sum_{l=0}^{m}\binom{m}{l}\frac{d^{m-l}}{dv^{m-l}}%
\Pr\left[  \left.  Y\left(  t+v\right)  =j\right\vert Y\left(  v\right)
=n\right]  \left(  \mathfrak{s}^{\left(  l\right)  }\left(  v\right)  \right)
_{n} \label{m-th_D_law_total_probability}%
\end{equation}
If we define the $N\times N$ matrix $\Phi\left(  t,v;k\right)  $ with elements
$\Phi_{jn}\left(  t,v;k\right)  =\frac{d^{k}}{du^{k}}\Pr\left[  \left.
Y\left(  t+v\right)  =j\right\vert Y\left(  v\right)  =n\right]  $, then
(\ref{m-th_D_law_total_probability}) is written in matrix form as%
\begin{equation}
\mathfrak{s}^{\left(  m\right)  }\left(  t+v\right)  =\sum_{l=0}^{m}\binom
{m}{l}\Phi\left(  t,v;m-l\right)  \mathfrak{s}^{\left(  l\right)  }\left(
v\right)  \label{m-th_D_law_total_probability_matrix}%
\end{equation}
Thus, the $m$-th derivative $\mathfrak{s}^{\left(  m\right)  }\left(
t+v\right)  $ of the probability state vector $\mathfrak{s}\left(  t+v\right)
$ at time $t+v$ is a linear combination of the vectors $\Phi\left(
t,v;m-l\right)  \mathfrak{s}^{\left(  l\right)  }\left(  v\right)  $,
illustrating that all the derivatives of the probability state vector
$\mathfrak{s}\left(  v\right)  $ at time $v$ are involved. In particular,
since $\Phi\left(  t,v;0\right)  =P\left(  t,v\right)  $, the transition
probability matrix with elements $\Pr\left[  \left.  Y\left(  t+v\right)
=j\right\vert Y\left(  v\right)  =n\right]  $ and $P\left(  t,0\right)  =$
$P\left(  t\right)  $ corresponds to the transition probability matrix of a
stationary continuous-time Markov process (i.e. independent of the time $v$).
Left-multiplying (\ref{m-th_D_law_total_probability_matrix}) with the vector
$u^{T}$ and invoking (\ref{D_axiom_1_prob}) shows that%
\[
u^{T}\Phi\left(  t,v;k\right)  =\left\{
\begin{array}
[c]{cc}%
u^{T} & \text{if }k=0\\
0 & \text{if }k>0
\end{array}
\right.
\]

Let us define the $mN\times1$ vector $\sigma\left(  t\right)  $ consisting of
the $m$-concatenation of the $N\times1$ vectors $\mathfrak{s}^{\left(
l\right)  }\left(  t\right)  $ for $0\leq l<m$, which is the column block
vector%
\begin{equation}
\sigma\left(  t\right)  =\left[
\begin{array}
[c]{c}%
\mathfrak{s}\left(  t\right) \\
\mathfrak{s}^{\left(  1\right)  }\left(  t\right) \\
\vdots\\
\mathfrak{s}^{\left(  m-1\right)  }\left(  t\right)
\end{array}
\right]  \label{def_sigma_vector}%
\end{equation}
then the vector relation (\ref{m-th_D_law_total_probability_matrix}) can be
written in matrix form as%
\begin{equation}
\sigma\left(  t+v\right)  =\digamma(t,v)\sigma\left(  v\right)
\label{General_stochastic_eq}%
\end{equation}
where the $mN\times mN$ square matrix
\[
\digamma(t,v)=\left[
\begin{array}
[c]{ccccc}%
\Phi\left(  t,v;0\right)  & O & O & \cdots & O\\
\Phi\left(  t,v;1\right)  & \Phi\left(  t,v;0\right)  & O & \cdots & O\\
\Phi\left(  t,v;2\right)  & \binom{2}{1}\Phi\left(  t,v;1\right)  &
\Phi\left(  t,v;0\right)  & \cdots & O\\
\vdots & \vdots & \vdots & \ddots & \\
\Phi\left(  t,v;m-1\right)  & \binom{m-1}{1}\Phi\left(  t,v;m-2\right)  &
\binom{m-1}{2}\Phi\left(  t,v;m-3\right)  & \cdots & \Phi\left(  t,v;0\right)
\end{array}
\right]
\]
is an lower triangular block matrix with determinant $\det\digamma
(t,v)=\left(  \det\Phi\left(  t,v;0\right)  \right)  ^{m}$ and whose inverse
is also a lower triangular block matrix as shown in
\cite{PVM_TriangularMatrix}. Thus,%
\begin{equation}
\sigma\left(  v\right)  =\digamma^{-1}(t,v)\sigma\left(  t+v\right)
\label{invers_general_stochatic}%
\end{equation}
If $t=0$ in (\ref{General_stochastic_eq}), then $\digamma(0,v)=I$ for any time
$t$. Thus, $\Phi\left(  0,v;0\right)  =I$ and $\Phi\left(  0,v;k\right)  =O$
for $k>0$.

\subsection{Matrix equation of fractional processes}

Returning to the general-order fractional, linear differential equation
(\ref{fractional_differential_state_vector_s(t)}) with general solution
$s_{\alpha}(t)=\sum_{l=0}^{m-1}t^{l}E_{\alpha,l+1}\left(  -Qt^{\alpha}\right)
s_{\alpha}^{\left(  l\right)  }\left(  0\right)  $ in
(\ref{general_solution_s_a(t)_Re(a)<m}), the $n$-th derivative is%
\[
s_{\alpha}^{\left(  n\right)  }(t)=\sum_{l=0}^{m-1}\frac{d^{n}}{dt^{n}}\left(
t^{l}E_{\alpha,l+1}\left(  -Qt^{\alpha}\right)  \right)  s_{\alpha}^{\left(
l\right)  }\left(  0\right)
\]
Invoking the key differential recursion \cite[art. 16]%
{PVM_Mittag-Leffler_Gamma} of the Mittag-Leffler function, $\frac{d^{n}%
}{dz^{n}}\left(  z^{b-1}E_{a,b}\left(  xz^{a}\right)  \right)  =z^{b-1-n}%
E_{a,b-n}\left(  xz^{a}\right)  $, yields%
\begin{equation}
s_{\alpha}^{\left(  n\right)  }(t)=\sum_{l=0}^{m-1}t^{l-n}E_{\alpha
,l-n+1}\left(  -Qt^{\alpha}\right)  s_{\alpha}^{\left(  l\right)  }\left(
0\right)  \label{n_derivative_state_vector_a}%
\end{equation}

Employing (\ref{n_derivative_state_vector_a}) for $n=m$ leads to the matrix
equation%
\begin{equation}
\sigma_{\alpha}\left(  t\right)  =\Upsilon(t,0)\sigma_{\alpha}\left(
0\right)  \label{general_matrix_solution_general-order_frac_lin_diff_eq}%
\end{equation}
where the block matrix
\[
\Upsilon(t,0)=\left[
\begin{array}
[c]{ccccc}%
E_{\alpha,1}\left(  -Qt^{\alpha}\right)  & tE_{\alpha,2}\left(  -Qt^{\alpha
}\right)  & t^{2}E_{\alpha,3}\left(  -Qt^{\alpha}\right)  & \cdots &
t^{m-1}E_{\alpha,m}\left(  -Qt^{\alpha}\right) \\
t^{-1}E_{\alpha,0}\left(  -Qt^{\alpha}\right)  & E_{\alpha,1}\left(
-Qt^{\alpha}\right)  & tE_{\alpha,2}\left(  -Qt^{\alpha}\right)  & \cdots &
t^{m-2}E_{\alpha,m-1}\left(  -Qt^{\alpha}\right) \\
t^{-2}E_{\alpha,-1}\left(  -Qt^{\alpha}\right)  & t^{-1}E_{\alpha,0}\left(
-Qt^{\alpha}\right)  & E_{\alpha,1}\left(  -Qt^{\alpha}\right)  & \cdots &
t^{m-3}E_{\alpha,m-2}\left(  -Qt^{\alpha}\right) \\
\vdots & \vdots & \vdots & \ddots & \\
t^{1-m}E_{\alpha,-m}\left(  -Qt^{\alpha}\right)  & t^{2-m}E_{\alpha
,1-m}\left(  -Qt^{\alpha}\right)  & t^{3-m}E_{\alpha,2-m}\left(  -Qt^{\alpha
}\right)  & \cdots & E_{\alpha,1}\left(  -Qt^{\alpha}\right)
\end{array}
\right]
\]
is not a lower triangular block matrix as $\digamma(t,0)$ in
(\ref{General_stochastic_eq}) for $v=0$.

\subsection{Matrix structure comparison between stochastic and fractional
process}

\label{sec_matrix_comparison}The general-order fractional, linear differential
equation (\ref{fractional_differential_state_vector_s(t)}) possesses a general
matrix solution $\sigma_{\alpha}\left(  t\right)  =\Upsilon(t,0)\sigma
_{\alpha}\left(  0\right)  $ in
(\ref{general_matrix_solution_general-order_frac_lin_diff_eq}) in the
generalized state vector $\sigma_{\alpha}\left(  t\right)  $ in
(\ref{def_sigma_vector}). If a general-order fractional, linear differential
equation (\ref{fractional_differential_state_vector_s(t)}) describes a
stochastic process on a finite state space, then its solution must also obey
the general matrix solution $\sigma\left(  t\right)  =\digamma(t,0)\sigma
\left(  0\right)  $ in (\ref{General_stochastic_eq}). Equating the generalized
state vectors $\sigma_{\alpha}\left(  t\right)  =\sigma\left(  t\right)  $
requires, for all times $t>0$, that the matrices $\Upsilon(t,0)$ and
$\digamma(t,0)$ have the same lower triangular block matrix structure.
However, the matrices $\Upsilon(t,0)$ and $\digamma(t,0)$ possess a different structure.

In summary, we must conclude that the general solution
(\ref{general_solution_s_a(t)_Re(a)<m}) does not describe a stochastic process
if $m=\left\lceil \alpha\right\rceil >1$.

\section{Conclusions}

\label{sec_conclusions}Our analysis indicates that the generalization of the
Markovian governing equation $\frac{ds\left(  t\right)  }{dt}=-Qs\left(
t\right)  $ in (\ref{differential_state_vector_s(t)}) by only replacing the
left-hand side $\frac{ds\left(  t\right)  }{dt}$ by a fractional derivative
$D_{m,0}^{\alpha}s\left(  t\right)  $, culminating in the generalized
governing equation (\ref{fractional_differential_state_vector_s(t)}), leads to
a mathematical incompatibility for $\alpha>1$. Since there is no obvious
reason why the fractional order $\alpha$ should be limited to the interval
$(0,1]$ and since measurements indicate that $\alpha>1$ for an overwhelming
amount of respiratory infectious diseases \cite{PVM_realdiseases_alpha>1}, we
believe that the formal process of replacing only the left-hand side of the
Markovian governing equation (\ref{differential_state_vector_s(t)}) with a
fractional derivative without modifying the right-hand side, leading to
$D_{m,0}^{\alpha}s_{\alpha}\left(  t\right)  =-Qs_{\alpha}\left(  t\right)  $,
needs a reconsideration. A mathematical equation acts as a balance and
changing only one-side logically destroys the balance. A consistent
generalization of the Markovian governing equation
(\ref{differential_state_vector_s(t)}) may need a modification of both sides
of the equation.

Since both equations (\ref{differential_state_vector_s(t)}) and
(\ref{fractional_differential_state_vector_s(t)}) possess the same right-hand
side $-Qs\left(  t\right)  $, their steady-state vector $s_{\infty}%
=\lim_{t\rightarrow\infty}s\left(  t\right)  =\lim_{t\rightarrow\infty
}s_{\alpha}\left(  t\right)  $ is the same, independent of the fractional
order $\alpha$ and satisfying $Qs_{\infty}=0$. Thus, if both processes start
(at $t=0$) with the same probability state vector $s\left(  0\right)
=s_{\alpha}\left(  0\right)  $, then they also end (at $t\rightarrow\infty$)
with the same probability state vector $s_{\infty}$, which poses a physical
constraint on the generalization of Markov processes.

Even the case for fractional order $\alpha\in(0,1]$, where
\emph{mathematically} no inconsistencies appear, is \emph{physically}
questionable, because the solution $s_{\alpha}\left(  t\right)  $ only seems
to be consistent with a semi-Markov process
\cite{PVM_fract_stoch_proc_spreading}. Thus, the case $\alpha\in(0,1]$ induces
an inconsistency: the fractional derivative
(\ref{general_in_m_complex_fractional_derivative_no_jump_alternative}) with
$p=0$ is a convolution of the process over all times since its start at $t=0$
and causes the current state at time $t$ to depend on its entire past history.
The semi-Markov process, on the other hand, only depends on the previous
states with Mittag-Leffler transition times. For epidemics on a contact graph,
a semi-Markov process causes two major problems (I) and (II), which contradict
physical observations. We quote from \cite{PVM_unraveling_frac_SIS}:

\textquotedblleft(I) all competing infection and curing processes are not
independent and (II) when a state transition occurs, all residual times of
competing infection and curing processes are neglected, due to the semi-Markov
\textquotedblleft restart\textquotedblright\ property. Dependence in (I)
implies that the curing process of some node $i$ in the contact graph directly
influences the infection and curing processes of every other node in the
network \cite{PVM_fracDerivative_Markov_epidemics},
\cite{PVM_fract_stoch_proc_spreading}. Restart (II) implies that, if a state
transition event occurs before node $i$ can cure, then the curing process of
node $i$ is reset as if that curing of node $i$ so far has never happened. The
beautiful memoryless property of the exponential and fundamental property in
Markov processes allows a reset, because the remaining time to cure is the
same as the time to cure after reset \cite[pp. 42-43, pp. 399-401]%
{PVM_PAComplexNetsCUP}. Both effects (I) and (II) are nonphysical in epidemic
spread.\textquotedblright

\section{Statements and Declarations}

\medskip AI software is not used at all. This research has been funded by the
European Research Council (ERC) under the European Union's Horizon 2020
research and innovation programme (grant agreement No 101019718).

{\footnotesize
\bibliographystyle{unsrt}
\bibliography{cac,math,misc,net,pvm,qth,tel}
}

\appendix{}

\section{Laplace transform of $\left.  \frac{d^{\alpha}f(z)}{dz^{\alpha}%
}\right\vert _{z=t}$ in
(\ref{general_in_m_complex_fractional_derivative_no_jump})}

\label{sec_Laplace_transform_complex_fractional_derivative}We summarize
results from \cite{PVM_fractional_derivative}. The single-sided\footnote{The
double-sided Laplace transform extends the lower integration bound of the
single-sided Laplace transform over all negative real values of $t$ and is
$\varphi(z)=\int_{-\infty}^{\infty}e^{-zt}f(t)dt$.} Laplace transform for
complex $z$ is defined (see e.g. \cite{Titchmarshfourier}, \cite[Chapter
VII]{Evgrafov_1965}, \cite{Widder})\ as
\begin{equation}
\varphi(z)=\int_{0}^{\infty}e^{-zt}f(t)dt \label{def_pgf_continuousrv}%
\end{equation}
and in operator form $\varphi(z)=\mathcal{L}_{z}\left[  f(t)\right]  $,
simplified to $\varphi(z)=\mathcal{L}\left[  f(t)\right]  $. The inverse
Laplace transform, in operator form $f(t)=\mathcal{L}_{t}^{-1}\left[
\varphi(z)\right]  $, simplified to $f(t)=\mathcal{L}^{-1}\left[
\varphi(z)\right]  $, is
\begin{equation}
f(t)=\frac{1}{2\pi i}\int_{c-i\infty}^{c+i\infty}\varphi(z)e^{zt}dz
\label{inverse_Laplace}%
\end{equation}
where $c$ is the smallest real value of $\operatorname{Re}(z)$ for which the
integral in (\ref{def_pgf_continuousrv}) converges.

The Laplace transform of the derivative\footnote{We assume that $f\left(
t\right)  $ is analytic at $t=0$, else $f\left(  0\right)  $ should be
replaced by $\lim_{t\downarrow0}f\left(  t\right)  =f\left(  0^{+}\right)  $,
because then the direction of the path towards the point $t=0$ matters and the
left-limit $\lim_{t\uparrow0}f\left(  t\right)  =f\left(  0^{-}\right)  $ and
right-limit may be different and the function $f\left(  t\right)  $ may not be
continuous at $t=0$.} $\mathcal{L}\left[  f^{\prime}(t)\right]  =\int
_{0}^{\infty}e^{-zt}f^{\prime}(t)dt=z\mathcal{L}\left[  f(t)\right]  -f\left(
0\right)  $. If we replace $f$ by $f^{\left(  k-1\right)  }$, then we find
with $y_{k}=\mathcal{L}\left[  f^{\left(  k\right)  }(t)\right]  $ the
recursion $y_{k}=zy_{k-1}-f^{\left(  k-1\right)  }\left(  0\right)  $, which
becomes after $l$ iterations%
\[
y_{k}=z^{l}y_{k-l}-\sum_{j=1}^{l}z^{j-1}f^{\left(  l-j\right)  }\left(
0\right)
\]
Since $y_{0}=\mathcal{L}\left[  f(t)\right]  $, choosing $l=k$ yields the
general form%
\[
\mathcal{L}\left[  f^{\left(  k\right)  }(t)\right]  =z^{k}\mathcal{L}\left[
f(t)\right]  -\sum_{j=1}^{k}z^{j-1}f^{\left(  k-j\right)  }\left(  0\right)
\]
The $k$-fold integral%
\begin{equation}
F_{k}\left(  t,a\right)  =\int_{a}^{t}du_{1}\int_{a}^{u_{1}}du_{2}\ldots
\int_{a}^{u_{k-1}}du_{k}f\left(  u_{k}\right)  \label{def_n_fold_integral}%
\end{equation}
can be reduced to a single integral, by repeated partial integration, for
$k\geq1$,%
\begin{equation}
F_{k}\left(  t,a\right)  =\frac{1}{\left(  k-1\right)  !}\int_{a}^{t}\left(
t-u\right)  ^{k-1}f\left(  u\right)  du \label{n_fold_integral}%
\end{equation}
The expression (\ref{n_fold_integral}), which possesses the form of
(\ref{general_in_m_complex_fractional_derivative_no_jump}), was known to
Cauchy, but is also contributed to Riemann and Liouville \cite[p.
43]{Duff_Naylor},\cite[Appendix E]{Gorenflo_2020}. Using the Riemann-Liouville
integral in (\ref{n_fold_integral}) for $k$-fold integration with $a=0$ and
the important convolution property $\mathcal{L}\left[  \left(  h\ast g\right)
\left(  t\right)  \right]  =\mathcal{L}\left[  \int_{0}^{t}h\left(
t-x\right)  g\left(  x\right)  dx\right]  =\mathcal{L}\left[  h\left(
t\right)  \right]  \mathcal{L}\left[  g\left(  t\right)  \right]  $ results
in
\[
\mathcal{L}\left[  F_{k}\left(  t,0\right)  \right]  =\mathcal{L}\left[
\frac{1}{\left(  k-1\right)  !}\int_{0}^{t}\left(  t-u\right)  ^{k-1}f\left(
u\right)  du\right]  =\frac{\mathcal{L}\left[  f(t)\right]  }{z^{k}}%
\]

Invoking the convolution property and $\mathcal{L}\left[  t^{-\beta}\right]
=\int_{0}^{\infty}t^{-\beta}e^{-zt}dt=\frac{\Gamma\left(  1-\beta\right)
}{z^{1-\beta}}$ valid for $\operatorname{Re}\left(  \beta\right)  <1$, the
Laplace transform (\ref{def_pgf_continuousrv}) of the generalization
(\ref{general_in_m_complex_fractional_derivative_no_jump}) with $p=0$ is%
\[
\mathcal{L}\left[  \left.  \frac{d^{\alpha}f(z)}{dz^{\alpha}}\right\vert
_{z=t}\right]  =\frac{1}{\Gamma\left(  m-\alpha\right)  }\mathcal{L}\left[
t^{-\alpha-1+m}\right]  \mathcal{L}\left[  f^{\left(  m\right)  }(t)\right]
\]
and%
\begin{equation}
\mathcal{L}\left[  \left.  \frac{d^{\alpha}f(z)}{dz^{\alpha}}\right\vert
_{z=t}\right]  =z^{\alpha}\left\{  \mathcal{L}\left[  f(t)\right]  -\sum
_{l=0}^{m-1}z^{-l-1}f^{\left(  l\right)  }\left(  0\right)  \right\}
\hspace{0.5cm}\text{for }\operatorname{Re}\left(  \alpha\right)  <m
\label{Laplace_transform_fractional_derivative_m_general}%
\end{equation}

\section{Proof of the law (\ref{Prob_t+u_in_u_tot_law_prob}) of total
probability for stochastic processes}

\label{sec_proof_law_total_prob}The definition of conditional probability,
$\Pr\left[  A\cap B\right]  =\Pr\left[  A|B\right]  \Pr\left[  B\right]  $,
illustrates that (\ref{Prob_t+u_in_u_tot_law_prob}) can be written as%
\begin{equation}
\Pr\left[  Y\left(  t+v\right)  =j\right]  =\sum_{n=1}^{N}\Pr\left[  \left\{
Y\left(  t+v\right)  =j\right\}  \cap\left\{  Y\left(  v\right)  =n\right\}
\right]  \label{Law_tot_prob_Joint_prob}%
\end{equation}
The case for $t=0$ is immediate, because $\left\{  Y\left(  v\right)
=j\right\}  \cap\left\{  Y\left(  v\right)  =n\right\}  =\emptyset$ for $j\neq
n$ as the process can only be in one state at time $v$. Now,%
\[
\Pr\left[  \cup_{n=1}^{N}\left\{  \left\{  Y\left(  t+v\right)  =j\right\}
\cap\left\{  Y\left(  v\right)  =n\right\}  \right\}  \right]  =\Pr\left[
\left\{  Y\left(  t+v\right)  =j\right\}  \cap\left\{  \cup_{n=1}^{N}Y\left(
v\right)  =n\right\}  \right]
\]
by the distributive property $\left\{  A\cap B_{1}\right\}  \cup\left\{  A\cap
B_{2}\right\}  =A\cap\left\{  B_{1}\cup B_{2}\right\}  $. Since $\cup
_{n=1}^{N}\left\{  Y\left(  v\right)  =n\right\}  =\Omega$ comprises the
entire sample space $\Omega$ and $A\cap\Omega=A$, we find that%
\[
\Pr\left[  Y\left(  t+v\right)  =j\right]  =\Pr\left[  \cup_{n=1}^{N}\left\{
\left\{  Y\left(  t+v\right)  =j\right\}  \cap\left\{  Y\left(  v\right)
=n\right\}  \right\}  \right]
\]
Since all events $W_{n}=\left\{  \left\{  Y\left(  t+v\right)  =j\right\}
\cap\left\{  Y\left(  v\right)  =n\right\}  \right\}  $ are mutually
exclusive\footnote{Associativity shows that%
\begin{align*}
W_{n}\cap W_{k}  &  =\left\{  \left\{  Y\left(  t+v\right)  =j\right\}
\cap\left\{  Y\left(  v\right)  =n\right\}  \right\}  \cap\left\{  \left\{
Y\left(  t+v\right)  =j\right\}  \cap\left\{  Y\left(  v\right)  =k\right\}
\right\} \\
&  =\left\{  Y\left(  t+v\right)  =j\right\}  \cap\left\{  \left\{  Y\left(
v\right)  =n\right\}  \cap\left\{  Y\left(  v\right)  =k\right\}  \right\}
=\left\{  Y\left(  t+v\right)  =j\right\}  \cap\emptyset=\emptyset
\end{align*}
} in that $W_{n}\cap W_{k}=\emptyset$ for $n\neq k$ and with relation
\begin{equation}
\Pr\left[  A\cup B\right]  =\Pr\left[  A\right]  +\Pr\left[  B\right]
-\Pr\left[  A\cap B\right]  \label{prob_union_events}%
\end{equation}
in \cite[eq. (2.4)]{PVM_PAComplexNetsCUP}, we arrive at
(\ref{Law_tot_prob_Joint_prob}) and thus prove
(\ref{Prob_t+u_in_u_tot_law_prob}).\hfill$\square\medskip$

\end{document}